\documentclass[10pt]{amsart}
\usepackage{amssymb}
\usepackage{amscd}
\usepackage{upref}

\theoremstyle{plain}
\newtheorem{thm}{Theorem}[section]
\newtheorem{cor}[thm]{Corollary}

\newtheorem{prop}[thm]{Proposition}

\newtheorem*{Theorem A}{Theorem A}
\newtheorem*{Theorem B}{Theorem B}
\newtheorem*{Theorem C}{Theorem C}

\theoremstyle{definition}

\theoremstyle{remark}

\numberwithin{equation}{section}

\begin{document}

\title[Infinite Product Decomposition of Orbifold Mapping Spaces]
{Infinite Product Decomposition of Orbifold Mapping Spaces}
\author{Hirotaka Tamanoi} 
\address[] {Department of Mathematics,
University of California Santa Cruz \newline 
\indent Santa Cruz, CA 95064} 
\email[]{tamanoi@math.ucsc.edu} 
\date{} 

\begin{abstract}
Physicists showed that the generating function of orbifold elliptic
genera of symmetric orbifolds can be written as an infinite
product. We show that there exists a geometric factorization on space
level behind this infinite product formula, and we do this in a much more general
framework of orbifold mapping spaces, where factors in the infinite product correspond to finite connected coverings of domain spaces whose fundamental groups are not necessarily abelian. From this formula, a concept of geometric Hecke operators for functors emerges. This is a non-abelian geometric generalization of usual Hecke operators. We show that these generalized Hecke operators indeed satisfy the identity of usual Hecke operators for the case of
2-dimensional tori.
\end{abstract}
\maketitle

\tableofcontents

\section{Introduction and summary of results}

 The elliptic genus of a Spin manifold $M$ refers to
the signature of $LM$ \cite{Land}, \cite{T3}. The elliptic genus of a complex manifold $M$ refers to the $S^1$-equivariant $\chi_y$-characteristic of its free loop space $LM=\text{Map}(S^1,M)$ \cite{H}. These are some of the versions of elliptic genera of $M$. Since $LM$ is infinite dimensional, the above statements must be be interpreted using a localization formula \cite{Wi}.

Let $G$ be a finite group. For a $G$-manifold $M$, we can consider an
orbifold version of the elliptic genus. However, the free loop space
$L(M/G)$ on the orbit space is not well behaved. Following \cite{HH}, we define 
the orbifold loop space $L_{\text{orb}}(M/G)$ by
\begin{equation}
L_{\text{orb}}(M/G)\overset{\text{def}}{=}
\bigl(\coprod_{g\in G}L_gM\bigr)/G
=\!\!\coprod_{(g)\in G_*}\!\!\bigl[L_gM/C_{G}(g)\bigr],
\label{1.1}
\end{equation}
where $G_*$ is the set of conjugacy classes in $G$, $C_G(g)$ is the
centralizer of $g$ in $G$, and $L_gM$ is the space of $g$-twisted
loops in $M$ given by 
\begin{equation}
L_gM=\{\gamma:\mathbb{R} \rightarrow M \mid \gamma(t+1)=g^{-1}\gamma(t) 
\text{ for all } t\in\mathbb{R}\}.
\label{1.2}
\end{equation}
The centralizer $C(g)$ acts on $L_gM$. Also note that if the order of
$g$ is finite and is equal to $s$, then each twisted
loop $\gamma$ in $L_gM$ is in fact a closed loop of length $s$. Thus,
$L_gM$ also admits an action of a circle $S^1=\mathbb{R}/s\mathbb{Z}$ of
length $s$.

One could use more sophisticated languages on orbifolds (see for example, \cite{M}), but for our purpose, the above definition suffices. 

Now the orbifold elliptic genus of $(M,G)$, denoted by
$\text{ell}_{\text{orb}}(M/G)$, is defined as the $S^1$-equivariant
$\chi_y$-characteristic of $L_{\text{orb}}(M/G)$:
\begin{equation}
\text{ell}_{\text{orb}}(M/G)
=\chi_y^{S^1}\bigl(L_{\text{orb}}(M/G)\bigr)
=\!\!\sum_{(g)\in G_*}\!\!\chi_y(L_gM)^{C(g)},
\label{1.3}
\end{equation}
where $\chi_y(L_gM)$ is thought of as
$R\bigl(C(g)\bigr)$-valued $S^1$-equivariant $\chi_y$-characteristic
computed and made sense through a use of localization formulae. 
Counting the dimension of coefficient vector spaces, we have 
\begin{equation}
\text{ell}_{\text{orb}}(M/G)\in\mathbb{Z}[y,y^{-1}][[q]],
\label{1.4}
\end{equation}
where the powers of $q$ are characters of $S^1$. 

Dijkgraaf, Moore, Verlinde and Verlinde \cite{DMVV} essentially proved a remarkable formula for the
generating function of orbifold elliptic genera of symmetric products.
This was subsequently extended to symmetric orbifold case by
Borisov-Libgober \cite{BL}. Here, for an integer $n\ge0$, the $n$-th
symmetric product of a space $X$ is defined as $SP^n(X)=X^n/\mathfrak{S}_n$,
where the $n$-th symmetric group $\mathfrak{S}_n$ acts on $X^n$ by
permuting $n$ factors. The DMVV and BL formula for the generating
function of orbifold elliptic genera of symmetric orbifolds is 
given by 
\begin{equation}
\begin{gathered}
\sum_{n\ge0}p^n\text{ell}_{\text{orb}}\bigl(SP^n(M/G)\bigr)
=\prod_{\substack{ n\ge1 \\m\ge0 \\ k\in\mathbb{Z} }}
(1-p^nq^my^k)^{-c(mn,k)},
\\
\text{where}\quad 
\text{ell}_{\text{orb}}(M/G)
=\sum_{\substack{ m\ge0 \\k\in\mathbb{Z}}}
c(m,k)q^my^k\in\mathbb{Z}[y,y^{-1}][[q]]. \label{1.5}
\end{gathered}
\end{equation}
The amazing thing about this formula is that the right hand side of
(1.5) is a genus 2 Siegel modular form, up to a simple multiplicative
factor. The main motivation of this paper is to understand a
geometric origin of this infinite product formula. In fact, we will
prove such an infinite product formula on a geometric level, not
merely on an algebraic level, as in (1.5).

We can describe this geometric formula in a general context. Let
$(M,G)$ be as before, and let $\Sigma$ be an arbitrary connected
manifold with $\Gamma=\pi_1(\Sigma)$. Instead of a loop space, we
consider a mapping space $\text{Map}(\Sigma,M/G)$. As before, this
space is not well behaved and the correct space to consider is the
orbifold mapping space defined by
\begin{equation}
\text{Map}_{\text{orb}}(\Sigma,M/G)
\overset{\text{def}}{=}
\Bigl(\!\!\!\!\!\!\!\!\!\coprod_{\theta\in\text{Hom}(\Gamma, G)}
\!\!\!\!\!\!\!\!\!\text{Map}_{\theta}(\widetilde{\Sigma},M)\Bigr)\Big/G
=\!\!\!\!\!\!\!\!\!\!\!\!\!\coprod_{(\theta)\in\text{Hom}(\Gamma,G)/G}
\!\!\!\!\!\!\!\!\!\!\!\!\!
\bigl[\text{Map}_{\theta}(\widetilde{\Sigma},M)/C(\theta)\bigr].
\label{1.6}
\end{equation}
Here $\widetilde{\Sigma}$ is the universal cover of $\Sigma$, and
$\text{Map}_{\theta}(\widetilde{\Sigma},M)$ is the space of
$\theta$-equivariant maps $\alpha:\widetilde{\Sigma} \rightarrow M$
such that $\alpha(p\cdot\gamma)=\theta(\gamma)^{-1}\cdot\alpha(p)$ for
all $p\in\widetilde{\Sigma}$ and $\gamma\in \Gamma$. Note here that we
regard the universal cover $\widetilde{\Sigma}$ as a
$\Gamma$-principal bundle over $\Sigma$. 

For a variable $t$ and a space $X$, let
$S_t(X)=\coprod_{k\ge0}t^kSP^k(X)$ be the total symmetric product of
$X$. For convenience, we often write this using the summation symbol
as $S_t(X)=\sum_{k\ge0}t^kSP^k(X)$. In this paper, summation symbol
applied to topological spaces means topological disjoint union. 

\begin{Theorem A}[Infinite Product Decomposition of Orbifold
Mapping Spaces of Symmetric Products] Let $M$ be a $G$-manifold and
let $\Sigma$ be a connected manifold. Then,
\begin{equation}
\sum_{n\ge0}p^n\text{\rm Map}_{\text{\rm orb}}\bigl(\Sigma,
SP^n(M/G)\bigr) \cong
\!\!\!\!\!\!\!\!\!\!\!\!\!\!
\prod_{\ \ \ \ [\Sigma' \to \Sigma]_{\text{\rm conn.}}}
\!\!\!\!\!\!\!\!\!\!\!\!\!\!
S_{p^{|\Sigma'/\Sigma|}}\bigl(\text{\rm Map}_{\text{\rm orb}}(\Sigma',
M/G)/\mathcal{D}(\Sigma'/\Sigma)\bigr).
\label{1.7}
\end{equation}
Here the infinite product is taken over all the isomorphism classes of
finite sheeted connected covering spaces $\Sigma'$ of $\Sigma$, and
$\mathcal{D}(\Sigma'/\Sigma)$ is the group of all deck transformations of
the covering space $\Sigma' \rightarrow\Sigma$ \textup{(}which is not
necessarily Galois\textup{)}. The number of sheets of this covering is
denoted by $|\Sigma'/\Sigma|$.
\end{Theorem A}
We will explain the details of the action of $\mathcal{D}(\Sigma'/\Sigma)$
on $\text{\rm Map}_{\text{\rm orb}}(\Sigma', M/G)$ in \S2. 

When $\Sigma=S^1$, the above formula reduces to 
\begin{equation}
\sum_{n\ge0}p^nL_{\text{orb}}\bigl(SP^n(M/G)\bigr)\cong
\prod_{r\ge1}S_{p^r}\bigl(L_{\text{orb}}^{(r)}(M/G)/\mathbb{Z}_r\bigr),
\label{1.8}
\end{equation}
where $L_{\text{orb}}^{(r)}(M/G)$ is the space of orbifold loops of
length $r$. This is the geometric version of the formula (1.5). This formula itself is relatively easy to prove. See \cite{WZ}. 

The above formula \eqref{1.8} for orbifold loop space is an "abelian" case since $\pi_1(S^1)\cong\mathbb Z$. The formula in Theorem A is, in a sense, a non-abelian generalization of this orbifold loop space case. The most interesting case seems to be the one in which $\Sigma$ is a $2$-dimensional surface (regarding it as a world-sheet of a moving string). Here, the genus of the surface can be arbitrary. In physics literature, elliptic genus itself is computed as a path integral over mapping spaces from torus \cite{DMVV}. 

Restricting the global decomposition formula \eqref{1.7} to the subspace of
constant orbifold maps and considering their numerical invariants, we
recover our previous results in \cite{T1,T2}. See section 3 for a
description of these results. We remark that we can apply
(generalized) homology and cohomology functors to \eqref{1.7} to obtain
infinite product decomposition formulas of these homology and
cohomology theories.

Another surprising formula discovered by physicists \cite{DMVV} is its
connection to Hecke operators. They showed that the right hand side of
formula \eqref{1.5} can be written in terms of Hecke operators in a very
nice way: 
\begin{equation}
\sum_{n\ge0}p^n\text{ell}_{\text{orb}}\bigl(SP^n(M/G)\bigr)
=\exp\Bigl(\sum_{r\ge1}p^rT(r)
\bigl[\text{ell}_{\text{orb}}(M/G)\bigr]\Bigr),
\label{1.9}
\end{equation}
where $T(r)$ is the $r$-th Hecke operator acting on weight 0 Jacobi
forms: 
\begin{equation}
T(r)\Bigl[\sum_{\substack{ m\ge0 \\ k\in\mathbb{Z}}}
c(m,k)q^my^{k}\Bigr]
=\sum_{ad=r}\frac 1a\sum_{\substack{ m\ge0 \\ k\in\mathbb{Z}}}
c(md,k)q^{am}y^{ak}.
\label{1.10}
\end{equation}
Is there a corresponding Hecke operator in our geometric context? Such
a Hecke operator must assign a certain space to a given space. Our
geometric decomposition formula \eqref{1.7} suggests what geometric Hecke
operators should be. For each positive integer $r$, we expect the
$r$-th {\it geometric Hecke operator} $\mathbb{T}(r)$ would act on a space
of the form $\text{Map}_{\text{orb}}(\Sigma, M/G)$, and produces a
space involving all the connected $r$-sheeted covering spaces of
$\Sigma$, as follows. 
\begin{equation}
\mathbb{T}(r)\bigl[\text{Map}_{\text{orb}}(\Sigma, M/G)\bigr]
\overset{\text{def}}{=}\!\!\!\!\!\!\!\!
\coprod_{\substack{ [\Sigma' \to \Sigma] \\ \ \ |\Sigma'/\Sigma|=r}}
\!\!\!\!\!\!\!\!
\text{Map}_{\text{orb}}(\Sigma', M/G)/\mathcal{D}(\Sigma'/\Sigma).
\label{1.11}
\end{equation}

The usual Hecke operators use covering spaces of the torus \cite{L}, and in \cite{DMVV}, they explain the above result \eqref{1.9} from this point of view. Our formula \eqref{1.11} uses covering spaces of $\Sigma$ whose fundamental group is not necessarily abelian. Thus, in a sense, our Hecke operator can be thought of as a non-abelian generalization of the usual Hecke operators.

A general discussion of geometric Hecke operators in the framework of
functors is more convenient and will be given in section 4.  Let
$\mathcal{F}$ be a functor from the category $\mathcal{C}$ of
topological spaces and continuous maps to itself. For example, for a
$G$-manifold $M$, let $\mathcal{F}_{(M,G)}$ be a conrtavariant functor
from $\mathcal{C}$ to itself given by $\mathcal{F}_{(M,G)}(\Sigma)
=\text{Map}_{\text{orb}}(\Sigma, M/G)$. Then, $\mathbb{T}(n)$ acts on
the functor $\mathcal{F}$ by the following formula for a connected
space $\Sigma$.
\begin{equation}
\bigl(\mathbb{T}(n)\mathcal{F}\bigr)(\Sigma)
\overset{\text{def}}{=}
\!\!\!\!\!\!\!\!\!\!\!\!
\coprod_{\substack{ \ \ \ \ [\Sigma' \rightarrow \Sigma]_{\text{conn.}} \\ 
|\Sigma'/\Sigma|=n }}
\!\!\!\!\!\!\!\!\!\!\!\!
\mathcal{F}(\Sigma')/\mathcal{D}(\Sigma'/\Sigma),
\label{1.12}
\end{equation}
where disjoint union runs over all isomorphism classes of connected
$n$-sheeted covering space of $\Sigma$. When $\Sigma$ is not
connected, we apply the above construction for each connected
component of $\Sigma$.  In terms of geometric Hecke operators, formula
\eqref{1.7} can be simply rewritten as
\begin{equation}
\sum_{n\ge0}p^n\text{Map}_{\text{orb}}\bigl(\Sigma, SP^n(M/G)\bigr)
\cong \prod_{r\ge1}S_{p^r}
\Bigl[\bigl(\mathbb{T}(r)\mathcal{F}_{(M,G)}\bigr)(\Sigma)\Bigr].
\tag{\ref{1.7}$'$}
\label{1.7'}
\end{equation}
It is very suggestive to compare this formula with \eqref{1.9}. If we
regard the $n$-th symmetric product $SP^n(X)$ as $X^n/n!$, since
$\mathfrak{S}_n$ has $n!$ elements, then we can regard $S_p(X)$ as
$\exp(pX)$. From this point of view, the analogy between \eqref{1.9} and
(1.$7'$) is reasonably precise. However, see also a remark after (4.2). 

The name geometric Hecke operator seems appropriate since these
operators do satisfy the usual identity when $\Sigma$ is a genus 1
Riemann surface.

\begin{Theorem B}[Hecke Identity for Geometric Hecke
Operators] Let $T$ be a $2$-dimensional torus. 
Let $\mathcal{F}$ be a functor from the category
$\mathcal{C}$ of topological spaces to itself. Then the geometric
Hecke operators $\mathbb{T}(n)$, $n\ge1$, satisfy
\begin{equation}
\bigl((\mathbb{T}(m)\circ\mathbb{T}(n))\mathcal{F}\bigr)(T)
=\!\!\!\sum_{d|(m,n)}\!\!\!
d\cdot \bigl(\bigl(\mathbb{T}\Big(\frac{mn}{d^2}\Bigr)\circ
\mathbb{R}(d)\bigr)\mathcal{F}\bigr)(T),
\label{1.13}
\end{equation}
where the operator $\mathbb{R}(d)$ on the functor $\mathcal{F}$ is given by 
\begin{equation}
(\mathbb{R}(d)\mathcal{F})(T)
=\mathcal{F}\bigl(R(d)T\bigr)
\big/\mathcal{D}\bigl(R(d)T/T\bigr),
\label{1.14}
\end{equation}
in which $R(d)T=\widetilde T/(d\cdot L)$ if $T=\widetilde{T}/L$ for some
lattice $L\subset \widetilde{T}\cong\mathbb{R}^2$. 
\end{Theorem B}

Thus, $R(d)T$ is a $d^2$-sheeted covering space of $T$. The
coefficient $d$ in the right hand side of \eqref{1.13} means a disjoint
topological union of $d$ copies. 

Note that \eqref{1.13} can be restated in a more familiar form as follows:
\begin{equation}
\begin{aligned}
\mathbb{T}(m)\circ\mathbb{T}(n)&=\mathbb{T}(mn), &\qquad&\text{if } (m,n)=1,\\
\mathbb{T}(p^r)\circ\mathbb{T}(p)
&=\mathbb{T}(p^{r+1})+p\cdot\mathbb{T}(p^{r-1})\circ\mathbb{R}(p),
&\ \ &\text{if $p$ prime}. 
\end{aligned}
\tag{\ref{1.13}$'$}\label{1.13'}
\end{equation}
As is well known in the theory of modular forms, these identities are
equivalent to an Euler product decomposition of the Dirichlet series with
the above Hecke operator coefficients. See \eqref{4.14}. 

It would be of interest to investigate relations among $\mathbb{T}(n)$s
when $\Sigma$ is a higher genus Riemann surfaces, or higher
dimensional tori whose fundamental group is free abelian. 

For a generalization of orbifold elliptic genus to the setting of generalized cohomology theory, see a paper by Ganter \cite{Ga}. 

The organization of this paper is as follows. In section 2, we prove
our main geometric decomposition formula in Theorem A. In section 3,
we specialize our infinite dimensional geometric formula to the finite
dimensional subspace of constant orbifold maps, and we deduce various
formulae of generating functions of orbifold invariants. In
section 4, after discussing some generality of geometric Hecke
operators on functors, we prove the Hecke identity \eqref{1.13}.

The main result of this paper, Theorem A, was first announced at a workshop at Banff International Research Station in June 2003.

\section{Infinite product
decomposition of orbifold mapping spaces}

First, we discuss some general facts of orbifold mapping spaces. For a
homomorphism $\theta:\Gamma \rightarrow G$ and a $\theta$-equivariant
map $\alpha:\widetilde{\Sigma} \rightarrow M$, let $\overline{\alpha}:
\Sigma \rightarrow M/G$ be the induced map on quotient spaces. Thus we
have a canonical map $\text{Map}_{\theta}(\widetilde{\Sigma}, M)
\rightarrow \text{Map}(\Sigma, M/G)$. Let $C_G(\theta)$ be the
centralizer of the image of $\theta$ in $G$. Note that inverse images of this map are $C_{G}(\theta)$ spaces. The action of $g\in G$ on $M$ has the effect
\begin{equation*}
g\cdot: \text{Map}_{\theta}(\widetilde{\Sigma}, M)
\longrightarrow 
\text{Map}_{g\cdot\theta\cdot g^{-1}}(\widetilde{\Sigma}, M),
\end{equation*}
and for every $\alpha\in\text{Map}_{\theta}(\widetilde{\Sigma},M)$, we
have $\overline{\alpha}=\overline{g\cdot\alpha}$ in
$\text{Map}(\Sigma, M/G)$.  Thus, we have a canonical map
\begin{equation}
\text{Map}_{\text{orb}}(\Sigma, M/G)\overset{\text{def}}{=}
\!\!\!\!\!\!\!\!\!\!\!\!\!\!\coprod_{(\theta)\in\text{Hom}(\Gamma,G)/G}
\!\!\!\!\!\!\!\!\!\!\!\!\!\!
\text{Map}_{\theta}(\widetilde{\Sigma}, M)/C_G(\theta) 
\longrightarrow \text{Map}(\Sigma, M/G).
\label{2.1}
\end{equation}
This map is in general not surjective nor injective. 

We consider a necessary condition for a map $f:\Sigma \rightarrow M/G$ to
have a lift to a $\theta$-equivariant map $\tilde{f}:
\widetilde{\Sigma} \rightarrow M$ for some $\theta$. Let $\eta$ be an
arbitrary contractible loop in $\Sigma$. Since $\widetilde{\Sigma}
\rightarrow\Sigma$ is a covering, $\eta$ always lifts to a contractible loop
$\tilde{\eta}$ in $\widetilde{\Sigma}$, and hence
$\tilde{f}(\tilde{\eta})$ is also contractible. Thus, for the
existence of a lift $\tilde{f}$ of a given map $f$, it is necessary
that for every contractible loop $\eta$ in $\Sigma$, $f(\eta)\subset
M/G$ lifts to a contractible loop in $M$.

Next, we discuss a functorial property of orbifold mapping spaces. 

\begin{prop} \textup{(i)} Let $M$ be a $G$-manifold. Any
map $f:\Sigma_1 \rightarrow \Sigma_2$ between connected manifolds induces a
well-defined map
\begin{equation}
f^*: \text{\rm Map}_{\textup{orb}}(\Sigma_2, M/G) \longrightarrow
\text{\rm Map}_{\textup{orb}}(\Sigma_1, M/G).  
\label{2.2} 
\end{equation}
\textup{(ii)} For two maps $f_1: \Sigma_1 \rightarrow \Sigma_2$ and 
$f_2: \Sigma_2 \rightarrow \Sigma_3$, we have $(f_2\circ f_1)^*=f_1^*\circ
f_2^*$. 
\end{prop}
\begin{proof} Let $\Gamma_i$ be the group 
$\mathcal{D}(\widetilde{\Sigma_i}/\Sigma_i)$ of all deck transformations
for the universal cover $\widetilde{\Sigma}_i \rightarrow \Sigma_i$ for
$i=1,2$. Since an isomorphism
$\mathcal{D}(\widetilde{\Sigma_i}/\Sigma_i) \cong \pi_1(\Sigma_i)$
depends on the choice of a base point in $\widetilde{\Sigma}_i$, it is
better to regard $\Gamma_i$ as the group of deck transformations
rather than as the fundamental group of $\Sigma_i$. We choose a lift
$\tilde{f}:\widetilde{\Sigma}_1 \rightarrow \widetilde{\Sigma}_2$ of
$f$. Then $\tilde{f}$ induces a homomorphism $\tilde{f}_*:\Gamma_1
\rightarrow \Gamma_2$ such that
$\tilde{f}(p\cdot\gamma_1)=\tilde{f}(p)\cdot\tilde{f}_*(\gamma_1)$ for
all $p\in \widetilde{\Sigma}_1$ and $\gamma_1\in\Gamma_1$. For a map
$\alpha\in \text{Map}_{\theta}(\widetilde{\Sigma}_2, M)$ with
$\theta\in \text{Hom}(\Gamma_2,G)$, we have $\alpha\circ \tilde{f}
\in\text{Map}_{\theta\circ\tilde{f}_*}(\widetilde{\Sigma}_1,
M)$. Hence the composition with $\tilde{f}$ gives an induced map
\begin{equation}
\tilde{f}^* : \!\!\!\!\!\!\!\!\!\!
\coprod_{\theta\in\text{Hom}(\Gamma_2,G)}
\!\!\!\!\!\!\!\!\!\!
\text{Map}_{\theta}(\widetilde{\Sigma}_2, M) \rightarrow 
\!\!\!\!\!\!\!\!\!\!\coprod_{\rho\in\text{Hom}(\Gamma_1,G)}
\!\!\!\!\!\!\!\!\!\!
\text{Map}_{\rho}(\widetilde{\Sigma}_1, M). 
\label{2.3}
\end{equation}
Obviously, this map commutes with the $G$-action on $M$. Hence by
quotienting by $G$, we have a map 
\begin{equation}
\tilde{f}^*: \text{Map}_{\text{orb}}(\Sigma_2, M/G) \rightarrow 
\text{Map}_{\text{orb}}(\Sigma_1, M/G).
\label{2.4}
\end{equation}
We have to verify that this map is independent of the chosen lift
$\tilde{f}$. Let $\tilde{f}':\widetilde{\Sigma}_1 \rightarrow
\widetilde{\Sigma}_2$ be another lift of $f$. By examining the image
of one point and using the uniqueness of lifts, we must have that
$\tilde{f}'=\tilde{f}\cdot \gamma_2$, globally on
$\widetilde{\Sigma}_1$, for some uniquely determined
$\gamma_2\in\Gamma_2$. Then, $(\alpha\circ\tilde{f}')(p_1)
=\alpha\bigl(\tilde{f}(p_1)\cdot\gamma_2\bigr)=
\theta(\gamma_2)^{-1}\cdot(\alpha\circ\tilde{f})(p_1)$
for all $p_1\in\widetilde{\Sigma}_1$. Note that $\theta(\gamma_2)\in
G$. Thus for all possible choices of lifts $\tilde{f}$, the collection
$\{\alpha\circ\tilde{f}\}$ is contained in a single $G$-orbit in 
$\coprod_{\rho\in\text{Hom}(\Gamma_1,G)}
\text{Map}_{\rho}(\widetilde{\Sigma}_1, M)$. Thus difference of
$\tilde{f}^*$ and $(\tilde{f}')^*$ in \eqref{2.3} disappear after dividing
by $G$, and the map \eqref{2.4} is independent of the choice of lifts
$\tilde{f}$. Hence we may simply call it $f^*$ as in \eqref{2.2}. 

The proof of the formula for the induced map of a composition is
routine. 
\end{proof}

As an immediate consequence, we have

\begin{cor} Let $\Sigma' \rightarrow \Sigma$ be a connected
covering space. Then the group $\mathcal{D}(\Sigma'/\Sigma)$ of all deck
transformations acts on $\text{\rm Map}_{\textup{orb}}(\Sigma', M/G)$. 
\end{cor}

For later use, we give details of this action. As before, let 
$\mathcal{D}(\widetilde{\Sigma}/\Sigma)=\Gamma$ and
$\Sigma'=\widetilde{\Sigma}/H$ for some $H\subset\Gamma$. Then $\mathcal{D}
(\Sigma'/\Sigma)\cong N_{\Gamma}(H)/H$. For
$f\in\text{Map}_{\rho}(\widetilde{\Sigma}', M)$, $u\in
N_{\Gamma}(H)$, and $g\in G$, the action of $u,g$ on $f$ is given by 
\begin{equation}
(u\cdot f)(p)=f(pu),\qquad (g\cdot f)(p)=g\cdot f(p), \quad
p\in\widetilde{\Sigma}'.
\label{2.5}
\end{equation}
These actions commute, but they do not preserve
$\rho\in\text{Hom}(H,G)$. How $\rho$ transforms under these actions
can be easily computed and we have the following commutative diagram:
\begin{equation}
\begin{CD}
\text{Map}_{\rho}(\widetilde{\Sigma}', M) @>{u\cdot}>{\cong}>
\text{Map}_{\rho^{u^{-1}}}(\widetilde{\Sigma}', M) \\
@V{g\cdot}V{\cong}V      @V{g\cdot}V{\cong}V \\
\text{Map}_{g\cdot \rho\cdot g^{-1}}(\widetilde{\Sigma}', M) 
@>{u\cdot}>{\cong}> 
\text{Map}_{g\cdot\rho^{u^{-1}}\cdot g^{-1}}(\widetilde{\Sigma}', M),
\end{CD}
\label{2.6}
\end{equation}
where $\rho^{u^{-1}}(h)=\rho(u^{-1}hu)$ for all $h\in H$. Since
$C_G(\rho)=C_G(\rho^{u^{-1}})$, commutativity of this diagram also
implies that for $u\in N_{\Gamma}(H)$, 
\begin{equation}
u\cdot: \text{Map}_{\rho}(\widetilde{\Sigma}', M) \xrightarrow{\cong} 
\text{Map}_{\rho^{u^{-1}}}(\widetilde{\Sigma}', M),\quad
\text{$C_G(\rho)$-equivariant}.
\label{2.7}
\end{equation}
A global statement is the following for $u\in N_{\Gamma}(H)$: 
\begin{equation}
u\cdot:\!\!\!\!\!\!\!\!\!\!\coprod_{\rho\in\text{Hom}(H,G)}
\!\!\!\!\!\!\!\!\!\!
\text{Map}_{\rho}(\widetilde{\Sigma}', M) \xrightarrow{\cong} 
\!\!\!\!\!\!\!\!\!\!
\coprod_{\rho\in\text{Hom}(H,G)}
\!\!\!\!\!\!\!\!\!\!
\text{Map}_{\rho}(\widetilde{\Sigma}', M),\quad
\text{$G$-equivariant}.
\label{2.8}
\end{equation}
In other words, the group $N_{\Gamma}(H)\times G$ acts on
$\coprod_{\rho}\text{Map}_{\rho}(\widetilde{\Sigma}', M)$. Also note
that the same group $N_{\Gamma}(H)\times G$ acts on the set
$\text{Hom}(H,G)$ by $[(u,g)\cdot\rho](h)=g\cdot\rho^{u^{-1}}(h)\cdot
g^{-1}$ for $h\in H$. The 
effect of changing $u\in N_{\Gamma}(H)$ by $h\in H$ can be computed as
\begin{equation}
\begin{aligned}
\rho^{(uh)^{-1}}(\ \cdot\ )
&=\rho(h)^{-1}\rho^{u^{-1}}(\ \cdot\ )\rho(h), \\
\rho^{(hu)^{-1}}(\ \cdot\ )
&=\rho^{u^{-1}}(h)^{-1}\rho^{u^{-1}}(\ \cdot\ )\rho^{u^{-1}}(h).
\end{aligned}
\label{2.9}
\end{equation}
This shows that modification of $u$ by elements in $H$ has the same
effect as the conjugation action by elements in $G$. Hence the map
induced from \eqref{2.8} on $G$-orbits is well defined for $\overline{u}\in
N_{\Gamma}(H)/H$, and we have
\begin{equation}
\overline{u}\cdot: \text{Map}_{\text{orb}}(\Sigma', M/G) \xrightarrow{\cong}
\text{Map}_{\text{orb}}(\Sigma', M/G).
\label{2.10}
\end{equation}
This is the action in Corollary 2.2. 

Since the action of $\mathcal{D}(\Sigma'/\Sigma)$ commutes with the
projection map $\pi:\Sigma' \rightarrow \Sigma$, the action of $\mathcal{D}
(\Sigma'/\Sigma)$ on $\text{Map}_{\text{orb}}(\Sigma', M/G)$ commutes
with the induced map $\pi^*$. In particular, the image of $\pi^*$ is in
the $\mathcal{D}(\Sigma'/\Sigma)$-fixed point subset: 
\begin{equation}
\text{Map}_{\text{orb}}(\Sigma, M/G) \xrightarrow{\pi^*}
\text{Map}_{\text{orb}}(\Sigma', M/G)^{\mathcal{D}(\Sigma'/\Sigma)}.
\label{2.11}
\end{equation}

We will need an identity on nested equivariant mapping spaces. Let $P
\rightarrow Z$ be a left $\Gamma$-equivariant right $G$-principal
bundle over a left $\Gamma$-space $Z$, where the left $\Gamma$-action
and the right $G$-action on $P$ commute. We simply call such a bundle
$\Gamma$-$G$ bundle \cite{May}. We studies this concept in detail
in section 3 of \cite{T2}, where the classification theorem of such
bundles is discussed. Note that $\text{Map}_G(P, M)$ is a left
$\Gamma$-space when $P$ is a $\Gamma$-$G$ bundle.

\begin{prop} With notations as above, we have 
\begin{equation}
\text{\rm Map}_{\Gamma}\bigl(\widetilde{\Sigma}, 
\text{\rm Map}_{G}(P, M)\bigr)
=\text{\rm Map}_{G}(\widetilde{\Sigma}
\underset{\Gamma}{\times} P, M).
\label{2.12}
\end{equation}
\end{prop}
\begin{proof} Without equivariance, this identity is obvious. So all we 
have to check is that the canonical correspondence preserves the
correct equivariance property.

Let $f:\widetilde{\Sigma} \rightarrow \text{Map}_G(P,M)$, and let
$u\in\widetilde{\Sigma}$. The $\Gamma$-equivariance of $f$ and
$G$-equivariance of $f(u)$ means $f(u\gamma)=\gamma^{-1}\cdot f(u)
=f(u)\circ\gamma$ and $f(u)(pg)=g^{-1}f(u)(p)$ for all
$\gamma\in\Gamma$, $g\in G$, $p\in P$. Let the canonically
corresponding map $\hat{f}:\widetilde{\Sigma}\times P \rightarrow M$ be
defined by $\hat{f}(u,p)=f(u)(p)$. The $\Gamma$-equivariance of $f$
implies that $\hat{f}(u\gamma,p)=\hat{f}(u,\gamma\cdot p)$ for all
$u,\gamma, p$. Hence $\hat{f}$ factors through
$\widetilde{\Sigma}\times_{\Gamma}P$ whose elements we denote by
$[u,p]$. Using $G$-equivariance of $f$, we have 
$\hat{f}([u,p]g)=\hat{f}([u,pg])=f(u)(pg)
=g^{-1}\cdot f(u)(p)=g^{-1}\hat{f}([u,p])$. Thus, $\hat{f}$ is
$G$-equivariant. 

The obvious inverse correspondence can be similarly checked to behave
correctly with respect to equivariance. 
\end{proof}

We examine the left hand side of the formula \eqref{1.7}. For a positive
integer $n$, let $\mathbf{n}=\{1,2,\dots, n\}$. Then the wreath product
$G_n=G\wr\mathfrak{S}_n$ is defined by 
\begin{equation}
G_n=G\wr\mathfrak{S}_n=\text{Map}(\mathbf{n}, G)\rtimes\mathfrak{S}_n.
\label{2.13}
\end{equation}
When $M$ is a $G$-manifold, the wreath product $G_n$ naturally acts on
the Cartesian product $M^n$, and its quotient space
$M^n/G_n=SP^n(M/G)$ is the $n$-the symmetric orbifold of $M/G$. For
detailed information on wreath product, see section 3 of \cite{T2}. To
understand \eqref{1.7}, first we note that
\begin{equation}
\text{Map}_{\text{orb}}\bigl(\Sigma, SP^n(M/G)\bigr)
=\!\!\!\!\!\!\!\!\!\!\!\!\!\!\!\!\!\!\!\!\!\!\!
\coprod_{\ \ \ \ \ \ \ (\theta)\in\text{Hom}(\Gamma,G_n)/G_n}
\!\!\!\!\!\!\!\!\!\!\!\!\!\!\!\!\!\!\!\!\!\!\!
\bigl[\text{Map}_{\theta}(\widetilde{\Sigma},M^n)
/C_{G_n}(\theta)\bigr].
\label{2.14}
\end{equation}
Let $\mathbf{n}\times G \rightarrow \mathbf{n}$ be the trivial $G$-principal
bundle over an $n$-element set $\mathbf{n}$. Since
$\text{Aut}_G(\mathbf{n}\times G)\cong G_n$ (see \cite{T2} Lemma 3-3),
the space of $G$-equivariant maps $\text{Map}_G(\mathbf{n}\times G,M)$
has the structure of left $G_n$ space and we have a $G_n$-equivariant
homeomorphism
\begin{equation}
M^n\cong \text{Map}_G(\mathbf{n}\times G,M).
\label{2.15}
\end{equation}
For a given homomorphism $\theta: \Gamma \rightarrow G_n$, both of the above
spaces can be thought of as $\Gamma$-spaces. Especially, the trivial
$G$-bundle $\mathbf{n}\times G \rightarrow \mathbf{n}$ acquires the structure of
a $\Gamma$-equivariant $G$-principal bundle, or simply a 
$\Gamma$-$G$ bundle, via $\theta$. We denote this by $(\mathbf{n}\times
G)_{\theta}$. Now \eqref{2.15} and Proposition 2.3 imply that
\begin{equation}
\text{Map}_{\theta}(\widetilde{\Sigma}, M^n)
\cong\text{Map}_{\Gamma}(\widetilde{\Sigma}, \text{Map}_G
\bigl((\mathbf{n}\times G)_{\theta},M)\bigr)
=\text{Map}_G\bigl(\widetilde{\Sigma}\times_{\Gamma}
(\mathbf{n}\times G)_{\theta},M\bigr).
\label{2.16}
\end{equation}

A $\Gamma$-$G$ bundle $P \rightarrow Z$ is called irreducible if $Z$ is a
transitive $\Gamma$-set. In this case, $\Gamma\times G$ acts
transitively on $P$. In section 3 of \cite{T2}, we classified all the
isomorphism classes of irreducible $\Gamma$-$G$ bundles. We showed
that any irreducible $\Gamma$-$G$ bundle must be of the form
$P_{H,\rho}=\Gamma\times_{\rho}G \rightarrow \Gamma/H$ for some subgroup
$H\subset\Gamma$ and a homomorphism $\rho:H \rightarrow G$. We also showed
that two irreducible $\Gamma$-$G$ bundles corresponding to
$(H_1,\rho_1)$ and $(H_2,\rho_2)$ are isomorphic as $\Gamma$-$G$
bundles if and only if (i) the subgroups $H_1$ and $H_2$ are conjugate
in $\Gamma$, and (ii) when $H_1=H_2=H$, we must have 
$[\rho_1]=[\rho_2]\in\text{Hom}(H,G)/(N_{\Gamma}(H)\times G)$ 
(\cite{T2}, Theorem E), where $N_{\Gamma}(H)$ and $G$ act on
$\text{Hom}(H,G)$ by conjugating $H$ and $G$, respectively. 

From now on, an element in $\text{Hom}(H,G)/(N_{\Gamma}(H)\times G)$ is
denoted with a square bracket as in $[\rho]$, and an element in
$\text{Hom}(H,G)/G$ is denoted by a round bracket as in $(\rho)$, to
distinguish these two kinds of conjugacy classes. 

Let $r_{\theta}(H,\rho)$ be the number of irreducible $\Gamma$-$G$
bundles isomorphic to $P_{H,\rho} \rightarrow \Gamma/H$ in the irreducible
decomposition of $(\mathbf{n}\times G)_{\theta} \rightarrow \mathbf{n}$. Thus, 
\begin{equation}
[(\mathbf{n}\times G)_{\theta} \rightarrow \mathbf{n}]
\cong\coprod_{[H]}\coprod_{[\rho]}\!\!\!\!
\coprod^{r_{\theta}(H,\rho)}\!\!\!\!
[P_{H,\rho} \rightarrow \Gamma/H].
\label{2.17}
\end{equation}
Here $[H]$ runs over all the conjugacy classes of finite index
subgroups of $\Gamma$, and for each $H$, $[\rho]$ runs over the set
$\text{Hom}(H,G)/(N_{\Gamma}(H)\times G)$. By examining the
decomposition of the base space $\mathbf{n}$ into transitive
$\Gamma$-sets, we have
\begin{equation}
\sum_{[H],[\rho]}r_{\theta}(H,\rho)|\Gamma/H|=n.
\label{2.18}
\end{equation}
Let $\mathbb{P}_{H,\rho}=\widetilde{\Sigma}\times_{\Gamma}P_{H,\rho}$ and
$\Sigma_H=\widetilde{\Sigma}\times_{\Gamma}(\Gamma/H)
=\widetilde{\Sigma}/H$. Then $\mathbb{P}_{H,\rho}$ is a $G$-bundle over a
covering space $\Sigma_H$ of $\Sigma$. Note that in $\mathbb{P}_{H,\rho}
\rightarrow \Sigma_H \rightarrow \Sigma$, for each point in $\Sigma$, fibres of
these bundles give $P_{H,\rho} \rightarrow \Gamma/H$. The above decomposition
now implies
\begin{equation}
\widetilde{\Sigma}\times_{\Gamma}
[(\mathbf{n}\times G)_{\theta} \rightarrow \mathbf{n}]
\cong\coprod_{[H]}\coprod_{[\rho]}\!\!\!\!
\coprod^{r_{\theta}(H,\rho)}\!\!\!\!
[\mathbb{P}_{H,\rho} \rightarrow \Sigma_H].
\label{2.19}
\end{equation}
This isomorphism allows us to rewrite (2-16) as 
\begin{equation}
\text{Map}_{\theta}(\widetilde{\Sigma},M^n)
\cong\prod_{[H]}\prod_{[\rho]}\!\!\!\!
\prod^{r_{\theta}(H,\rho)}\!\!\!\!
\text{Map}_G(\mathbb{P}_{H,\rho}, M)
\cong\prod_{[H]}\prod_{[\rho]}\!\!\!\!
\prod^{r_{\theta}(H,\rho)}\!\!\!\!
\text{Map}_{\rho}(\widetilde{\Sigma}_H, M).
\label{2.20}
\end{equation}
The last isomorphism is because $\mathbb{P}_{H,\rho}
=\widetilde{\Sigma}_H\times_{\rho}G$.  This gives multiplicative
decomposition of each disjoint summand of the right hand side of
\eqref{2.14}. Next, we need to understand the centralizer
$C_{G_n}(\theta)$ of the image of the homomorphism $\theta: \Gamma
\rightarrow G_n$ in $G_n$. One of the main results of \cite{T2} is the
description of the structure of the centralizer $C_{G_n}(\theta)$. It
says that
\begin{equation}
C_{G_n}(\theta)\cong\prod_{[H]}\prod_{[\rho]}
\bigl[\text{Aut}_{\text{$\Gamma$-$G$}}(P_{H,\rho})
\wr\mathfrak{S}_{r_{\theta}(H,\rho)}\bigr],
\label{2.21}
\end{equation}
where $\text{Aut}_{\text{$\Gamma$-$G$}}(P_{H,\rho})$ is the group of
$\Gamma$-equivariant $G$-principal bundle automorphisms of $P_{H,\rho}
\rightarrow \Gamma/H$. In terms of the $G$-bundle $\mathbb{P}_{H,\rho} \rightarrow
\Sigma_H$ over a covering space,
$\text{Aut}_{\text{$\Gamma$-$G$}}(P_{H,\rho})$ is isomorphic to the
group $\text{Aut}_{G}(\mathbb{P}_{H,\rho})_{\Sigma_H/\Sigma}$ of $G$-bundle
isomorphisms of $\mathbb{P}_{H,\rho}$ whose induced map on $\Sigma_H$ is a
deck transformation of $\Sigma_H \rightarrow \Sigma$ (\cite{T2}, Proposition
7-3). 

Next we describe the structure of $\text{Aut}_{\text{$\Gamma$-$G$}}
(P_{H,\rho})$. We recall that the group $N_{\Gamma}(H)\times G$ acts
on the set $\text{Hom}(H,G)$ by $(u,g)\cdot \rho=
g\cdot\rho^{u^{-1}}\cdot g^{-1}$ for $u\in N_{\Gamma}(H)$, $g\in G$ and
$\rho\in\text{Hom}(H,G)$. Let $T_{\rho}$ be the isotropy subgroup of
this action at $\rho$:
\begin{equation}
T_{\rho}=\{(u,g)\in N_{\Gamma}(H)\times G \mid
g\cdot\rho^{u^{-1}}(h)\cdot g^{-1}=\rho(h) 
\text{ for all }h\in H \}.
\label{2.22}
\end{equation}
This group $T_{\rho}$ contains a subgroup $H_{\rho}
=\bigl\{\bigl(h,\rho(h)\bigr)\in T_{\rho} \mid h\in H\bigr\}\cong
H$. Then Theorem 4-4 in \cite{T2} shows that $H_{\rho}$ is a normal
subgroup of $T_{\rho}$ and we have the following exact sequence:
\begin{equation}
1 \rightarrow H_{\rho} \rightarrow T_{\rho} \rightarrow 
\text{Aut}_{\text{$\Gamma$-$G$}}
(P_{H,\rho}) \rightarrow 1.
\label{2.23}
\end{equation}
Now we are ready to prove Theorem A. 

\begin{proof}[Proof of Theorem A]
Using \eqref{2.14}, \eqref{2.18}, \eqref{2.20}, \eqref{2.21}, 
we can rewrite the left hand side of \eqref{1.8} as 
\begin{align*}
\sum_{n\ge0}
&p^n\text{Map}_{\text{orb}}\bigl(\Sigma, SP^n(M/G)\bigr)\\
&=\sum_{n\ge0}\sum_{[\theta]}\prod_{[H]}\prod_{[\rho]}
p^{r_{\theta}(H,\rho)|\Gamma/H|}\Bigl[\bigl(
\!\!\!\!\!\!\!
\prod^{\ \ \ r_{\theta}(H,\rho)}
\!\!\!\!\!\!\!\!\text{Map}_{\rho}
(\widetilde{\Sigma}_H, M)\bigr)
/\bigl(\text{Aut}_{\text{$\Gamma$-$G$}}(P_{H,\rho})
\wr\mathfrak{S}_{r_{\theta}(H,\rho)}\bigr)\Bigr] \\
\intertext{Here $\text{Aut}_{\text{$\Gamma$-$G$}}(P_{H,\rho})
\cong \text{Aut}_{G}(\mathbb{P}_{H,\rho})_{\Sigma_H/\Sigma}$ acts on 
$\text{Map}_{\rho}(\widetilde{\Sigma}_H, M)\cong 
\text{Map}_G(\mathbb{P}_{H,\rho}, M)$ by the obvious action.}
&=\sum_{n\ge0}\sum_{[\theta]}\prod_{[H]}\prod_{[\rho]}
p^{r_{\theta}(H,\rho)|\Gamma/H|}SP^{r_{\theta}(H,\rho)}
\bigl(\text{Map}_{\rho}(\widetilde{\Sigma}_H, M)
/\text{Aut}_{\text{$\Gamma$-$G$}}(P_{H,\rho})\bigr)\\
&=\prod_{[H]}\prod_{[\rho]}\Bigl[\sum_{r\ge0}p^{r|\Gamma/H|}SP^r
\bigl(\text{Map}_{\rho}(\widetilde{\Sigma}_H, M)
/\text{Aut}_{\text{$\Gamma$-$G$}}(P_{H,\rho})\Bigr]\\
&=\prod_{[H]}\prod_{[\rho]}
S_{p^{|\Gamma/H|}}
\bigl(\text{Map}_{\rho}(\widetilde{\Sigma}_H, M)
/\text{Aut}_{\text{$\Gamma$-$G$}}(P_{H,\rho})\bigr)\\
&=\prod_{[H]}S_{p^{|\Gamma/H|}}
\Bigl[\coprod_{[\rho]}
\text{Map}_{\rho}(\widetilde{\Sigma}_H, M)
/\text{Aut}_{\text{$\Gamma$-$G$}}(P_{H,\rho})\Bigr].
\end{align*}
Here in the above formulae,
$[\rho]\in\text{Hom}(H,G)/(N_{\Gamma}(H)\times G)$. On the other hand, 
since $\mathcal{D}(\Sigma_H/\Sigma)\cong N_{\Gamma}(H)/H$, we have
\begin{align*}
\text{Map}_{\text{orb}}(\Sigma_H, M/G)
/\mathcal{D}(\Sigma_H/\Sigma)
&=\bigl[\bigl(\!\!\!\!\!\!\!\!\!\!\!\!\!\!\!
\coprod_{\ \ \ \ \ \rho\in\text{Hom}(H,G)}
\!\!\!\!\!\!\!\!\!\!\!\!\!\!\!\!
\text{Map}_{\rho}(\widetilde{\Sigma}_H, M)\bigr)\big/G\bigr]
\big/\bigl(N_{\Gamma}(H)/H\bigr) \\
&=\bigl(\!\!\!\!\!\!\!\!\!\!\!\!\!\!\!
\coprod_{\ \ \ \ \ \rho\in\text{Hom}(H,G)}
\!\!\!\!\!\!\!\!\!\!\!\!\!\!\!\!
\text{Map}_{\rho}(\widetilde{\Sigma}_H, M)\bigr)
/(N_{\Gamma}(H)\times G).
\end{align*}
Here we recall that the action of $G$ and $N_{\Gamma}(H)$ commutes,
and the action of $H\subset N_{\Gamma}(H)$ can be absorbed into the
action of $G$. See \eqref{2.5}, \eqref{2.6}, \eqref{2.8} and
\eqref{2.9} for details on this. In particular, the action of
$(u,g)\in N_{\Gamma}(H)\times G$ is such that
\begin{equation*}
(u,g): \text{Map}_{\rho}(\widetilde{\Sigma}_H, M) \xrightarrow{\cong}
\text{Map}_{g\rho^{u^{-1}}g^{-1}}(\widetilde{\Sigma}_H, M).
\end{equation*}
Since $T_{\rho}$ in \eqref{2.22} is exactly the subgroup which preserves
$\rho\in\text{Hom}(H,G)$ under $(N_{\Gamma}(H)\times G)$-action, in the
above identity, we get
\begin{equation*}
\text{Map}_{\text{orb}}(\Sigma_H, M/G)
/\mathcal{D}(\Sigma_H/\Sigma)
=\coprod_{[\rho]}\bigl(\text{Map}_{\rho}(\widetilde{\Sigma}_H, M)
/T_{\rho}\bigr),
\end{equation*}
where $[\rho]$ runs over the orbit set $\text{Hom}(H,G)
/(N_{\Gamma}(H)\times G)$. Next observe that the subgroup $H_{\rho}$
of $T_{\rho}$ acts trivially on
$\text{Map}_{\rho}(\widetilde{\Sigma}_H, M)$. To see this, let
$\bigl(h,\rho(h)\bigr)\in H_{\rho}$ for $h\in H$, and $f\in
\text{Map}_{\rho}(\widetilde{\Sigma}_H, M)$. Then, for any
$p\in\widetilde{\Sigma}_H$, we have 
\begin{equation*}
\bigl[\bigl(h,\rho(h)\bigr)f\bigr](p)
=\rho(h)\cdot(hf)(p)
=\rho(h)f(ph)=\rho(h)\rho(h)^{-1}f(p)=f(p). 
\end{equation*}
Thus, $H_{\rho}$ acts trivially on
$\text{Map}_{\rho}(\widetilde{\Sigma}_H, M)$. Hence quotienting by
$T_{\rho}$ in the above formula can be replaced by quotienting by
$T_{\rho}/H_{\rho}\cong
\text{Aut}_{\text{$\Gamma$-$G$}}(P_{H,\rho})$. Thus, collecting all the
above calculations, we finally have 
\begin{equation*}
\sum_{n\ge0}p^n\text{Map}_{\text{orb}}\bigl(\Sigma, SP^n(M/G)\bigr)
=\prod_{[H]}S_{p^{|\Gamma/H|}}
\bigl(\text{Map}_{\text{orb}}(\Sigma_H, M/G)
/\mathcal{D}(\Sigma_H/\Sigma)\bigr).
\end{equation*}
This completes the proof. 
\end{proof}

When $G=\{1\}$, we have $\text{Map}_{\text{orb}}(\Sigma,
M)=\text{Map}(\Sigma, M)$, and formula (1-8) becomes 
\begin{equation}
\sum_{n\ge0}p^n\text{\rm Map}_{\text{\rm orb}}\bigl(\Sigma,
SP^n(M)\bigr) \cong
\!\!\!\!\!\!\!\!\!\!\!\!\!\!
\prod_{\ \ \ \ [\Sigma' \to \Sigma]_{\text{\rm conn.}}}
\!\!\!\!\!\!\!\!\!\!\!\!\!\!
S_{p^{|\Sigma'/\Sigma|}}\bigl(\text{\rm Map}(\Sigma',
M)/\mathcal{D}(\Sigma'/\Sigma)\bigr).
\label{2.24}
\end{equation}

\section{Generating functions of finite orbifold invariants}

We specialize our main decomposition formula of infinite dimensional
orbifold mapping spaces to the finite dimensional subspace of constant
orbifold maps. Most of the results in \cite{T1, T2} follow from this
restricted formula, and we reproduce some of the main results in these
papers as corollaries to Theorem A.

Since $\text{Map}_{\theta}(\widetilde{\Sigma}, M)_{\text{const.}}\cong
M^{\langle\theta\rangle}$, where $M^{\langle\theta\rangle}$ denotes
the fixed point subset of $\theta$, we have
\begin{equation}
\text{Map}_{\text{orb}}(\Sigma, M/G)_{\text{const.}}
=\!\!\!\!\!\!\!\!\!\!\!\!\!\!\!\!\!\!
\coprod_{\ \ \ \ (\theta)\in\text{Hom}(\Gamma, G)/G}
\!\!\!\!\!\!\!\!\!\!\!\!\!\!\!\!\!\!
\bigl[M^{\langle\theta\rangle}/C(\theta)\bigr]
\overset{\text{def}}{=}
C_{\Gamma}(M/G).
\label{3.1}
\end{equation}
As an immediate consequence of Theorem A, we have the following
decomposition formula for constant orbifold maps. 

\begin{prop} Let $M$ be a $G$-space and let $\Gamma$
be an arbitrary group. Then, 
\begin{equation}
\begin{aligned}
\sum_{n\ge0}p^nC_{\Gamma}\bigl(SP^n(M/G)\bigr)
&=\prod_{[H]}S_{p^{|\Gamma/H|}}
\bigl(C_H(M/G)/(N_{\Gamma}(H)/H)\bigr)\\
&=\prod_{[H]}S_{p^{|\Gamma/H|}}
\bigl(\coprod_{[\rho]}(M^{\langle\rho\rangle}/T_{\rho})\bigr),
\end{aligned}
\label{3.2}
\end{equation}
where $[H]$ runs over all the conjugacy classes of finite index
subgroups of $\Gamma$, and for each $[H]$, $[\rho]$ runs over the set 
$\text{\rm Hom}(H,G)/\bigl(N_{\Gamma}(H)\times G\bigr)$. 
\end{prop}
Note that in Theorem A, $\Gamma$ is the fundamental group of the
manifold $\Sigma$. But after eliminating $\Sigma$ by considering
constant orbifold maps, $\Gamma$ can be an arbitrary (discrete) group
in Proposition 3.1.

Here we comment on the action of $N_{\Gamma}(H)/H$ on
$C_H(M/G)=(\coprod_{\rho}M^{\langle \rho\rangle})/G$ in \eqref{3.2}, where
$\rho\in\text{Hom}(H,G)$. In view of \eqref{2.5}, the action of
$N_{\Gamma}(H)$ commutes with the action of $G$, and for any $u\in
N_{\Gamma}(H)$ and any $x\in M^{\langle
\rho\rangle}$, the action of $u$ on $x$ is such that $u\cdot x=x$, as
can be easily verified by \eqref{2.5}. However, this does not mean that the
action of $N_{\Gamma}(H)$ on $C_H(M/G)$ is trivial. in fact, it is not
trivial in general. What happens is that the action of $u$ sends
$M^{\langle \rho\rangle}$ to $M^{\langle \rho^{u^{-1}}\rangle}$, where
$G$-conjugacy classes $(\rho)$ and $(\rho^{u^{-1}})$ can be distinct,
although these two spaces are identical subspaces of
$M$, since $\langle\rho\rangle=\langle\rho^{u^{-1}}\rangle$ as
subgroups of $G$. For a given $(\rho)\in\text{Hom}(H,G)/G$, let
$N_{\Gamma}^{\rho}(H)$ be the isotropy subgroup of $N_{\Gamma}(H)$ at
$(\rho)$. Recall that we have an exact sequence of groups [\cite{T2},
formula (4-6)]:
\begin{equation*}
1 \rightarrow C_G(\rho) \rightarrow T_{\rho} \rightarrow 
N_{\Gamma}^{\rho}(H) \rightarrow 1.
\end{equation*}
Thus,
$M^{\langle\rho\rangle}/T_{\rho}
=\bigl(M^{\langle\rho\rangle}/C(\rho)\bigr)/N_{\Gamma}^{\rho}(H)$. 
We examine the action of $u\in N_{\Gamma}^{\rho}(H)$ on
$M^{\langle
\rho\rangle}/C(\rho)$. By definition, for any $u\in
N_{\Gamma}^{\rho}(H)$, $\rho$ and $\rho^{u^{-1}}$ are $G$-conjugate, 
and thus there exists $g\in G$ such that
$\rho^{u^{-1}}(h)=g^{-1}\rho(h)g$ for all $h\in H$. This means that
$(u,g)\in T_{\rho}$. We have 
\begin{equation*}
M^{\langle \rho\rangle}/C(\rho) \xrightarrow{u\cdot=\text{Id}}
M^{\langle \rho^{u^{-1}}\rangle}/C(\rho^{u^{-1}})
=M^{\langle g^{-1}\rho g\rangle}/C( g^{-1}\rho g)
\xrightarrow[\cong]{g\cdot}
M^{\langle \rho\rangle}/C(\rho),
\end{equation*}
by \eqref{2.6}. This means that when we apply $u\cdot$, $\rho$ moves within
the same $G$-conjugacy class to $\rho^{u^{-1}}$. To bring it back to
$\rho$, we then apply $g\in G$. Thus, for $u\in N_{\Gamma}^{\rho}(H)$
and $x\in M^{\langle \rho\rangle}/C(\rho)$, the action of $u$ on $x$
is given by $u\cdot x=g\cdot x$ where $g\in G$ is an arbitrary element
such that $(u,g)\in T_{\rho}$.

Let $\chi(X)$ be the topological Euler characteristic for a
topological space $X$. In \cite{T2}, we introduced a notion of an
orbifold Euler characteristic associated to a group $\Gamma$ defined
for a $G$-manifold $M$:
\begin{equation}
\chi_{\Gamma}(M;G)
\overset{\text{def}}{=}
\chi\bigl(C_{\Gamma}(M/G)\bigr)
=\!\!\!\!\!\!\!\!\!\!\!\!\!\!\!\!\!\!
\sum_{\ \ \ \ (\theta)\in\text{Hom}(\Gamma,G)/G}
\!\!\!\!\!\!\!\!\!\!\!\!\!\!\!\!\!\!
\chi\bigl(M^{\langle\theta\rangle}/C(\theta)\bigr).
\label{3.3}
\end{equation}
We observe that when $\Gamma=\mathbb{Z}$, 
\begin{equation}
\chi_{\mathbb{Z}}(M;G)=\!\!\!\!\!\sum_{\ \ (g)\in G_*}\!\!\!\!\!
\chi\bigl(M^{\langle g\rangle}/C(g)\bigr)
=\frac1{|G|}\sum_{gh=hg}\chi\bigl(M^{\langle g,h\rangle}\bigr)
\label{3.4}
\end{equation}
is the physicist's orbifold Euler characteristic $e_{\text{orb}}(M/G)$
\cite{DHVW}. Here in the last summation, the pair $(g,h)$ runs over
the set of commuting pairs of elements. The second identity is due to
Lefschetz Fixed Point Formula. Formula \eqref{3.3} gives the correct
generalization of $e_{\text{orb}}(M/G)$ since it comes from a very
natural geometry of orbifold mapping spaces \eqref{3.1}.

In \cite{T2}, we introduced a notion of orbifold Euler
characteristic of $M/G$ associated to a $\Gamma$-set $X$, denoted by
$\chi_{[X]}(M;G)$. When $X$ is a transitive $\Gamma$-set of the form
$X=\Gamma/H$, it is given by
\begin{equation}
\begin{gathered}
\chi_{[\Gamma/H]}(M;G)
=\chi\bigl(C_H(M/G)/(N_{\Gamma}(H)/H)\bigr) \\
\begin{aligned}
\text{where}\quad
C_H(M/G)/(N_{\Gamma}(H)/H)
&=\!\!\!\!\!\!\!\!\!\!\!\!\!\!\!\!\!\!\!\!\!\!\!\!\!\!\!\!\!\!\!\!\!\!
\!\!\!\!
\coprod_{\ \ \ \ \ \ \ \ \ \ \ \ [\rho]\in\text{Hom}(H,G)
/(N_{\Gamma}(H)\times G)}
\!\!\!\!\!\!\!\!\!\!\!\!\!\!\!\!\!\!\!\!\!\!\!\!\!\!\!\!\!\!\!\!\!\!
\!\!\!\!
M^{\langle\rho\rangle}
/\text{Aut}_{\text{$\Gamma$-$G$}}(P_{H,\rho}) \\
&=\coprod_{[\rho]}M^{\langle\rho\rangle}/T_{\rho}.
\end{aligned}
\end{gathered}
\label{3.5}
\end{equation}
The second identity above can be proved on topological space level by
an argument similar to the last part of the proof of Theorem A. 

Now we compute the topological Euler characteristic of both sides of
\eqref{3.2}. We recall that $\chi\bigl(S_p(X)\bigr)=(1-p)^{-\chi(X)}$.

\begin{cor}[\cite{T2} Theorem C] Let $M$ be a $G$-set and
let $\Gamma$ be an arbitrary group. The the generating function of
orbifold Euler characteristic associated to $\Gamma$ of symmetric
orbifolds is given by
\begin{equation}
\sum_{n\ge0}p^n\chi_{\Gamma}(M^n;G_n)
=\prod_{[H]}(1-p^{|\Gamma/H|})
^{-\chi_{[\Gamma/H]}(M;G)},
\label{3.6}
\end{equation}
where $[H]$ runs over all conjugacy classes of finite index subgroups
of $\Gamma$. 
\end{cor} 
We can rewrite \eqref{3.6} in terms of Hecke operators as follows. For a
$G$-manifold, let $\chi_{(M;G)}$ be an integer valued function on the
set of discrete groups given by 
\begin{equation}
\chi_{(M;G)}(\Gamma)
\overset{\text{def}}{=}
\chi\bigl(C_{\Gamma}(M/G)\bigr)=\chi_{\Gamma}(M;G).
\label{3.7}
\end{equation}
For an integer $n\ge1$, let a Hecke operator $\mathbb{T}(n)$ act on the
function $\chi_{(M;G)}$ by 
\begin{equation}
\bigl[\mathbb{T}(n)\chi_{(M;G)}\bigr](\Gamma)
\overset{\text{def}}{=}
\!\!\!\!\!\!\!
\sum_{\substack{ [H] \\ \ \ |\Gamma/H|=n }}
\!\!\!\!\!\!\!
\chi\bigl(C_H(M/G)\big/(N_{\Gamma}(H)/H)\bigr),
\label{3.8}
\end{equation}
so that $\mathbb{T}(n)\chi_{(M;G)}$ is another integral function on the
set of discrete groups. Then as functions on the set of groups, \eqref{3.6}
means
\begin{equation}
\sum_{n\ge0}p^n\chi_{(M^n;G_n)}
=\prod_{n\ge1}(1-p^n)^{-\mathbb{T}(n)\chi_{(M;G)}}.
\label{3.9}
\end{equation}

Now we consider the case in which $\Gamma$ is {\it abelian}. In this
case, the action of $N_{\Gamma}(H)=\Gamma$ on $H\subset \Gamma$ is
trivial and so dividing by $N_{\Gamma}(H)/H$ has no effect. Thus, we
have $C_H(M/G)/(N_{\Gamma}(H)/H)=C_H(M/G)$ and
consequently, 

\begin{cor} Let $\Gamma$ be an arbitrary abelian
group. For any $G$-space $M$, we have 
\begin{equation}
\sum_{n\ge0}p^n\chi_{\Gamma}(M^n;G_n)
=\prod_{H}(1-p^{|\Gamma/H|})^{-\chi_H(M;G)},
\label{3.10}
\end{equation}
where the product is over all finite index subgroups $H$ of $\Gamma$. 
\end{cor}

In particular, when $\Gamma=\mathbb{Z}$, the formula \eqref{3.7} reduces to 
\begin{equation}
\sum_{n\ge0}p^ne_{\text{orb}}\bigl(SP^n(M/G)\bigr)
=\prod_{r\ge1}(1-p^r)^{-e_{\text{orb}}(M/G)}.
\label{3.11}
\end{equation}
This is the formula proven in \cite{HH} when $G$ is trivial, and for
general $G$ in \cite{W}. 

Instead of Euler characteristic, we can consider other numerical
invariants such as signature, spin index, $\chi_y$-characteristic,
etc., in suitable categories of manifolds. The formula \eqref{3.2} will then
provide us with infinite product formula of the corresponding
generating functions of orbifold invariants of symmetric
orbifolds. What is more interesting in this context is that, since we
have a decomposition on the space level, we can apply various
(generalized) homology and cohomology functors to obtain infinite
product decomposition formulae. This will be discussed in future
papers.

\section{Geometric Hecke operators for functors}

In this section, we prove the Hecke identity \eqref{1.13} for
2-dimensional tori.  Let $\mathcal{C}$ be the category of topological
spaces and continuous maps. Let $\mathcal{F}:\mathcal{C} \rightarrow
\mathcal{C}$ be a covariant (or contravariant) functor. Then it
formally follows that whenever $f:X \rightarrow Y$ is a homeomorphism,
the corresponding map $\mathcal{F}(f): \mathcal{F}(X) \rightarrow
\mathcal{F}(Y)$ (or $\mathcal{F}(Y) \rightarrow \mathcal{F}(X)$ in the
contravariant case) is also a homeomorphism. In particular, this
implies that when $X$ is a $G$-space, it automatically follows that
$\mathcal{F}(X)$ is also a $G$-space.

The geometric Hecke operator $\mathbb{T}(n)$, $n\ge1$, acts on 
a functor $\mathcal{F}$ as follows. For any connected space $X\in \mathcal{C}$,
\begin{equation}
\bigl(\mathbb{T}(n)\mathcal{F}\bigr)(X)\overset{\text{def}}{=}
\!\!\!\!\!\!\!\!\!\!\!\!\!\!
\coprod_{\substack{ \ \ \ [X' \rightarrow X]_{\text{conn.}} \\ |X'/X|=n}}
\!\!\!\!\!\!\!\!\!\!\!\!\!\!
\mathcal{F}(X')/\mathcal{D}(X'/X),
\label{4.1}
\end{equation}
where the disjoint union runs over the isomorphism classes of
connected $n$-sheeted covering spaces $X'$ of $X$, and $\mathcal{D}(X'/X)$
is the group of all deck transformations of $X' \rightarrow X$. When $X$ is
not connected, we apply the above construction to each of the
connected component. 

In general, we do not expect $\mathbb{T}(n)\mathcal{F} : \mathcal{C}
\rightarrow \mathcal{C}$ to be a functor. However, see Proposition 4.1
where such a situation does occur.

For the purpose of this paper, the main example of the functor
$\mathcal{F}$ is of course the orbifold mapping space functor. Namely,
for any $G$-space $M$, and any connected space $\Sigma$, we let
\begin{equation*}
\mathcal{F}_{(M;G)}(\Sigma)=\text{Map}_{\text{orb}}(\Sigma, M/G).
\end{equation*}
Proposition 2-1 shows that this is indeed a contravariant functor in
$\Sigma$. In terms of this notation, Theorem A can be restated as a
formal power series of functors as 
\begin{equation}
\sum_{n\ge0}p^n\mathcal{F}_{(M^n;G_n)}
=\prod_{n\ge1}S_{p^n}\bigl(\mathbb{T}(n)\mathcal{F}_{(M;G)}\bigr).
\label{4.2}
\end{equation}
However, in some context, for example in the Grothendieck ring of
varieties, it can make sense and can be justified to write
$S_p(X)=(1-p)^{-X}$ using powers whose exponents are spaces
\cite{GLM}. For the purpose of our present paper, we can regard
$S_p(X)$ as the definition of $(1-p)^{-X}$. This is more appropriate for our purpose since, for example, for Euler characteristic, we have $\chi\bigl(S_p(X(\bigr)=(1-p)^{-\chi(X)}$ for any space $X$. In this point of view, Theorem A has the following form:
\begin{equation}
\sum_{n\ge0}p^n\mathcal{F}_{(M^n;G_n)}(\Sigma)
=\prod_{n\ge1}(1-p^n)^{-(\mathbb{T}(n)\mathcal{F}_{(M;G)})(\Sigma)}.
\label{4.3}
\end{equation}
By Proposition 4.1 below, this formula can be regarded as a generating
function of functors from the category $\mathcal{C}_{\pi_1}$ to
$\mathcal{C}$, where $\mathcal{C}_{\pi_1}$ is the category of
topological spaces whose morphisms are restricted to those continuous
maps inducing isomorphisms on fundamental groups.

By considering constant orbifold maps, we have $\mathcal{F}_{(M;G)}
(\Sigma)_{\text{const.}}=C_{\pi_1(\Sigma)}(M/G)$. Then, by taking
topological Euler characteristic of \eqref{4.3} restricted to constant
orbifold maps, we recover the formula \eqref{3.9}. Notice that factors
$(1-p^n)$ in \eqref{3.9} are already present in \eqref{4.3} on space
level.

To define a composition of geometric Hecke operators, we need to have
functoriality of geometric Hecke operators in a certain special
situation. 

\begin{prop} Let $\mathcal{F}:\mathcal{C} \rightarrow \mathcal{C}$ be a
covariant functor. Let $X$ and $Y$ be connected spaces, and let $f: X
\rightarrow Y$ be a map such that $f_*:\pi_1(X) \rightarrow \pi_1(Y)$ is an
isomorphism. Then for every positive integer $n$, $f$ induces a map
\begin{equation}
f_*: (\mathbb{T}(n)\mathcal{F})(X) \rightarrow 
(\mathbb{T}(n)\mathcal{F})(Y), 
\label{4.4}
\end{equation}
such that for $ X \xrightarrow{f} Y \xrightarrow{g} Z$, we have
$(g\circ f)_*=g_*\circ f_*$.

A similar statement holds for contravariant functors. 
\end{prop}
\begin{proof} We fix a base point $x_0$ of $X$. Let 
$p: X' \rightarrow X$ be a connected $n$-sheeted covering space. For
each choice of a base point $x_0'$ of $X'$ over $x_0$, the subgroup
$H=p_*\bigl(\pi_1(X',x_0')\bigr)$ has index $n$ in
$\pi_1(X,x_0)$. Since $f_*: \pi_1(X,x_0) \rightarrow \pi_1(Y,y_0)$,
where $y_0=f(x_0)$, is an isomorphism by hypothesis, the subgroup
$f_*(H)$ has index $n$ in $\pi_1(Y,y_0)$. Let $(Y',y_0')$ be a
connected $n$-sheeted covering space with base point corresponding to
$f_*(H)$. The choice of $y_0'$ is unique up to the action of the group
$\mathcal{D}(Y'/Y)$ of deck transformations. Note that since $f_*:
\pi_1(X,x_0) \rightarrow \pi_1(Y,y_0)$ is an isomorphism, $f_*$
induces an isomorphism between the corresponding deck transformations
$\mathcal{D}(X'/X) \xrightarrow[\cong]{f_*} \mathcal{D}(Y'/Y)$. By the
Lifting Theorem in covering space theory, there exists a unique
$\mathcal{D}(X'/X)$-equivariant map $\tilde{f}: X' \rightarrow Y'$
such that $\tilde{f}(x_0')=y_0'$. By the functorial property, we see
that $\mathcal{F}(\tilde{f}): \mathcal{F}(X') \rightarrow \mathcal{F}(Y')$
is $\mathcal{D}(X'/X)\cong \mathcal{D}(Y'/Y)$-equivariant. Hence it induces
a map on the quotient:
\begin{equation*}
\overline{\mathcal{F}(\tilde{f})} : \mathcal{F}(X')/\mathcal{D}(X'/X)
\rightarrow \mathcal{F}(Y')/\mathcal{D}(Y'/Y).
\end{equation*}
Different choices of the lift $\tilde{f}$ are related by the action of
deck transformations. Hence the map $\overline{\mathcal{F}(\tilde{f})}$ on
the orbit space depends only on $f$. Repeating the above constructions
for each isomorphism class of connected $n$-sheeted covering spaces of
$X$, we obtain a map
\begin{equation}
f_*: \!\!\!\!\!\!\!\!\!\!\!\!\!\!
\coprod_{\substack{ \ \ \ [X' \rightarrow X]_{\text{conn.}} \\ |X'/X|=n}}
\!\!\!\!\!\!\!\!\!\!\!\!\!\!
\mathcal{F}(X')/\mathcal{D}(X'/X) \rightarrow 
\!\!\!\!\!\!\!\!\!\!\!\!\!\!
\coprod_{\substack{ \ \ \ [Y' \rightarrow Y]_{\text{conn.}} \\ |Y'/Y|=n }}
\!\!\!\!\!\!\!\!\!\!\!\!\!\!
\mathcal{F}(Y')/\mathcal{D}(Y'/Y).
\label{4.5}
\end{equation}
This is the map \eqref{4.4}. The behavior under the composition of two maps
can be easily verified. The argument for contravariant functors is
similar. 
\end{proof}

As a special case, let $f: X \rightarrow X$ be a homeomorphism. There
is one point which we have to be careful about in the above
construction of $f_*$. For a connected $n$-sheeted covering space
$p:(X',x_0') \rightarrow (X,x_0)$, the based covering space $(X'',y_0)
\rightarrow (X,f(x_0))$ corresponding to the subgroup
$f_*\bigl(p_*\bigl(\pi_1(X',x_0')\bigr)\bigr)\subset
\pi_1\bigl(X,f(x_0)\bigr)$ may not be isomorphic to $X' \rightarrow X$
as a covering space over $X$, although $X'$ and $X''$ are homeomorphic
via a lift $\tilde{f}: X' \xrightarrow{\cong} X''$ of $f$. Thus, in general,
the induced map
\begin{equation}
f_* :  \!\!\!\!\!\!\!\!\!\!\!\!\!\!
\coprod_{\substack{ \ \ \ [X' \rightarrow X]_{\text{conn.}} \\ |X'/X|=n}}
\!\!\!\!\!\!\!\!\!\!\!\!\!\!
\mathcal{F}(X')/\mathcal{D}(X'/X) \xrightarrow{\cong}
\!\!\!\!\!\!\!\!\!\!\!\!\!\!
\coprod_{\substack{ \ \ \ [X' \rightarrow X]_{\text{conn.}} \\ |X'/X|=n}}
\!\!\!\!\!\!\!\!\!\!\!\!\!\!
\mathcal{F}(X')/\mathcal{D}(X'/X)
\label{4.6}
\end{equation}
shuffles connected components, and it is not easy to control this
shuffling. This is an obstacle in studying compositions of Hecke
operators given in \eqref{4.7} below. However, when $f:X \rightarrow
X$ is a deck transformation of some covering $X \rightarrow X_0$, the
situation can be completely clarified. In particular, when
$\pi_1(X_0)$ is abelian, it turns out that the action of
$\mathcal{D}(X/X_0)$ on $(\mathbb{T}(n)\mathcal{F})(X)$ does preserve
connected components, and there is a simple relation among various
groups of deck transformations involved.

Anyway, as a formal consequence of Proposition 4.1, we have

\begin{cor} Let $\mathcal{F}:\mathcal{C} \rightarrow \mathcal{C}$ be an
arbitrary covariant or contravariant functor. If $X$ is $G$-space,
then for every positive integer $n$, the space
$(\mathbb{T}(n)\mathcal{F})(X)$ is also a $G$-space.  
\end{cor}

Next, we consider compositions of Hecke operators given as follows. 
\begin{equation}
\begin{aligned}
\bigl((\mathbb{T}(m)\circ \mathbb{T}(n))\mathcal{F}\bigr)(X)
&=\mathbb{T}(m)\bigl(\mathbb{T}(n)\mathcal{F}\bigr)(X) \\
&=\!\!\coprod_{\ [X']_m}\!\!
\bigl[\bigl(\mathbb{T}(n)\mathcal{F}\bigr)(X')\bigr]
\big/\mathcal{D}(X'/X) \\
&=\!\!\coprod_{\ [X']_m}\!\!
\Bigl[\!\!\!\coprod_{\ \ [X'']_n}\!\!\!
\mathcal{F}(X'')/\mathcal{D}(X''/X')\Bigr]\Big/\mathcal{D}(X'/X),
\end{aligned}
\label{4.7}
\end{equation}
where $[X']_m$ runs over the set of isomorphism classes of connected
$m$-sheeted covering spaces of $X$, and for a given $X'$, $[X'']_n$
runs over the set of isomorphism classes of connected $n$-sheeted
covering spaces of $X'$.  

As remarked earlier concerning formula \eqref{4.6}, the action of the
group of deck transformations $\mathcal{D}(X'/X)$ on
$\bigl(\mathbb{T}(n)\mathcal{F}\bigr)(X')$ permutes its connected
components. We now clarify what happens.

Let $\widetilde{X} \rightarrow X$ be the universal cover of $X$ and let
$\Gamma=\mathcal{D}(\widetilde{X}/X)\cong\pi_1(X)$ be its group of deck
transformations. We regard $\widetilde{X} \rightarrow X$ as the right
$\Gamma$-principal bundle over $X$. Let $K\subset H\subset \Gamma$ be
subgroups such that $|\Gamma/H|=m$ and $|H/K|=n$. We put
$X_K=\widetilde{X}/K$ and $X_H=\widetilde{X}/H$. Then $X_H \rightarrow X$ is
a connected $m$-sheeted covering of $X$ with $\mathcal{D}(X_H/X)\cong
N_{\Gamma}(H)/H$, and $X_K \rightarrow X_H$ is a connected $n$-sheeted
covering of $X_H$ with $\mathcal{D}(X_K/X_H)\cong N_{H}(K)/K$. Let $g\in
N_{\Gamma}(H)\subset\Gamma$. Then the right multiplication by $g$
induces the following diagram of homeomorphisms and covering spaces:
\begin{equation}
\begin{CD}
\widetilde{X} @>>> X_K @>>> X_H @>>> X \\
@V{\cdot g}V{\cong}V  @V{\cdot g}V{\cong}V @V{\cdot g}V{\cong}V @| \\
\widetilde{X} @>>> X_{g^{-1}Kg} @>>> X_H @>>> X.
\end{CD}
\label{4.8}
\end{equation} 
Since $g\in N_{\Gamma}(H)$, the map $\cdot g: X_H \xrightarrow{\cong}
X_H$ is a deck transformation of $X_H$ over $X$. However, since $g$
may not be in $N_{\Gamma}(K)$, $\cdot g: X_K \xrightarrow{\cong}
X_{g^{-1}Kg}$ is only an isomorphism of covering spaces over $X$. If
$g\in N_{\Gamma}(K)$, then $X_{g^{-1}Kg}=X_K$ and $\cdot g$ induces a
deck transformation of $X_K$ over $X$. For the middle square, when
$g\in H\subset N_{\Gamma}(H)$, $\cdot g$ induces an isomorphism of two
coverings $X_K$ and $X_{g^{-1}Kg}$ over $X_H$. If, furthermore, we
have $g\in N_{H}(K)\subset H$, then $\cdot g$ induces a deck
transformation of $X_K$ over $X_H$. This clarifies the action of
$\mathcal{D}(X'/X)$ on $(\mathbb{T}(n)\mathcal{F})(X')$ where $X'=X_H$ and
$X''=X_K$.

The above situation simplifies when the fundamental group of $X$ is
abelian. In this case, every element $g\in\Gamma$ induces a deck
transformation $\cdot g: X_H \xrightarrow{\cong} X_H$ whose lift
$\cdot g: X_K \xrightarrow{\cong} X_K$ preserves $X_K$. Also we have
$\mathcal{D}(X_K/X_H)\cong H/K$ for any two subgroups $K\subset
H\subset \Gamma$. The formula \eqref{4.7} now simplifies as follows.

\begin{prop} Let $X$ be a connected space whose
fundamental group is abelian. Then, the composition of two geometric
Hecke operators is given by 
\begin{equation}
\bigl(\mathbb{T}(m)(\mathbb{T}(n)\mathcal{F})\bigr)(X)
=\!\!\!\!\!
\coprod_{\substack{ H\subset\Gamma \\|\Gamma/H|=m }}
\Bigl[
\!\!\!\!\!\!
\coprod_{\substack{ K\subset H \\ \ \ \ |H/K|=n}}
\!\!\!\!\!\!\!\!\!
\bigl(\mathcal{F}(X_K)/\mathcal{D}(X_K/X)\bigr)\Bigr].
\label{4.9}
\end{equation}
\end{prop}
\begin{proof} By \eqref{4.7}, we have 
\begin{align*}
\bigl(\mathbb{T}(m)(\mathbb{T}(n)\mathcal{F})\bigr)(X)\bigr)
&=\!\!\!\!\!  \coprod_{\substack{ H\subset\Gamma \\|\Gamma/H|=m }}
\Bigl[ \!\!\!\!\!\!  \coprod_{\substack{ K\subset H \\ \ \ \ |H/K|=n
}} \!\!\!\!\!\!\!\!\!
\bigl(\mathcal{F}(X_K)/\mathcal{D}(X_K/X_H)\bigr)\Bigr]
\Big/\mathcal{D}(X_H/X)\\ 
\intertext{Since $\Gamma$ is abelian,
$\mathcal{D}(X_H/X)$ preserves $\mathcal{F}(X_K)/\mathcal{D}(X_K/X_H)$
for each $K\subset H$,} 
&=\!\!\!\!\!  \coprod_{\substack{ H\subset\Gamma
\\|\Gamma/H|=m }} \!\!\!\!\!\!  \coprod_{\substack{ K\subset H \\ \ \ \
|H/K|=n}} \!\!\!\!\!\!\!\!\!
\Bigl[\bigl(\mathcal{F}(X_K)/\mathcal{D}(X_K/X_H)\bigr)
\Big/\mathcal{D}(X_H/X)\Bigr] \\ \intertext{since
$\mathcal{D}(X_K/X_H)=H/K$, $\mathcal{D}(X_H/X)=\Gamma/H$, and
$\mathcal{D}(X_K/X)=\Gamma/K$, we have} &=\!\!\!\!\!  \coprod_{\substack{
H\subset\Gamma \\|\Gamma/H|=m }} \Bigl[ \!\!\!\!\!\!  \coprod_{\substack{
K\subset H \\ \ \ \ |H/K|=n}} \!\!\!\!\!\!\!\!\!
\bigl(\mathcal{F}(X_K)/\mathcal{D}(X_K/X)\bigr)\Bigr].  
\end{align*}
This completes the proof. 
\end{proof}

We continue to assume that the fundamental group of $X$ is
abelian. For an integer $d\ge1$, let $R(d)X$ be the covering space of
$X$ corresponding to $d\cdot\pi_1(X)\subset\pi_1(X)$. Let
$\mathbb{R}(d)$ act on a functor $\mathcal{F}$ by
\begin{equation}
(\mathbb{R}(d)\mathcal{F})(X)
\overset{\text{def}}{=}
\mathcal{F}\bigl(R(d)X\bigr)\big/\mathcal{D}\bigl((R(d)X)/X\bigr).
\label{4.10}
\end{equation}
As in Proposition 4.1, we can show that any map $f:X \rightarrow X$
inducing an isomorphism on fundamental groups gives rise to a map 
\begin{equation}
f_*:(\mathbb{R}(d)\mathcal{F})(X) \rightarrow 
(\mathbb{R}(d)\mathcal{F})(X).
\label{4.11}
\end{equation}
In particular, if $X$ is a $G$-space, then not only $\mathcal{F}(X)$ is a
$G$-space, but also $(\mathbb{R}(d)\mathcal{F})(X)\bigr)$ is a
$G$-space for all $d\ge1$. 

The main result in this section is the following Hecke identity for
geometric Hecke operators for 2-dimensional tori $T$. 

\begin{thm} Let $\mathcal{F}:\mathcal{C} \rightarrow \mathcal{C}$
be an arbitrary contravariant or covariant functor. Let $T$ be a
$2$-dimensional torus. Then for every pair of positive integers $m$
and $n$, the composition of two geometric Hecke operators satisfy
\begin{equation}
\bigl(\mathbb{T}(m)(\mathbb{T}(n)\mathcal{F})\bigr)(T)
=\!\!\!\!\sum_{d|(m,n)}\!\!\!\!d\cdot
\bigl(\mathbb{T}\left(\frac{mn}{d^2}\right)
(\mathbb{R}(d)\mathcal{F})\bigr)(T). 
\label{4.12}
\end{equation}
In particular, $\mathbb{T}(m)$ and $\mathbb{T}(n)$ commute. 
\end{thm}
In the right hand side of \eqref{4.9}, the summation symbol means disjoint
topological union, and the factor $d$ means a disjoint union of $d$
copies. 

For the proof, we first recall the ordinary Hecke identity for
lattices. For details, see (\cite{L}, p.16). Let $\mathcal{A}$ be the
free abelian group generated by rank $2$ lattices $L$ of $\mathbb{C}$. The
the Hecke operator $T(n)$ for $n\ge1$ is a map $T(n): \mathcal{A}
\rightarrow \mathcal{A}$ defined by
\begin{equation*}
T(n)(L)=\!\!\!\!\!\!\!\sum_{\ \ [L:L']=n}\!\!\!\!\!\!\!
L'\in\mathcal{A}.
\end{equation*}
Let $R(n):\mathcal{A} \rightarrow \mathcal{A}$ be defined by
$R(n)L=nL$ consisting of elements $\{n\cdot\ell\}_{\ell\in L}\subset
L$. Then Hecke identity says
\begin{equation}
T(m)\circ T(n)(L)=\!\!\!\!\sum_{d|(m,n)}\!\!\!\!
d\cdot R(d)\circ T\left(\frac{mn}{d^2}\right)(L),
\label{4.13}
\end{equation}
for any lattice $L$. From this formula, it is clear that $T(m)$ and
$T(n)$ commute. Also, it is easy to check that $R(d)$ and $T(n)$
commute. 

Now we are ready to prove Theorem 4.4.

\begin{proof}[Proof of Theorem 4.4] Since $\Gamma\cong\pi_1(T)\cong
\mathbb{Z}^2$ is free abelian of rank $2$, any subgroup of $\Gamma$ of
finite index is also free abelian of rank $2$. Applying \eqref{4.9} in
our context, we obtain
\begin{equation*}
\bigl(\mathbb{T}(m)(\mathbb{T}(n)\mathcal{F})\bigr)(T)
=\!\!\!\!\!
\coprod_{\substack{ H\subset\Gamma \\|\Gamma/H|=m }}
\Bigl[
\!\!\!\!\!\!
\coprod_{\substack{ K\subset H \\ \ \ \ |H/K|=n }}
\!\!\!\!\!\!\!\!\!
\bigl(\mathcal{F}(T_K)/(\Gamma/K)\bigr)\Bigr],
\end{equation*}
where $T_K$ is a covering torus corresponding to an index $mn$
sublattice $K\subset\Gamma$. By the ordinary Hecke identity \eqref{4.13},
any index $mn$ sublattice $K$ of $\Gamma$ arising in the above
disjoint union is of the form $d\cdot L$ for some integer $d$ dividing
$(m,n)$, and for some lattice $L$ of index $(mn)/d^2$ in $\Gamma$, and
furthermore, there are exactly $d$ such sublattices in the above
disjoint union. Hence the right hand side of the above expression can
be rewritten as
\begin{equation*}
\bigl(\mathbb{T}(m)(\mathbb{T}(n)\mathcal{F})\bigr)(T)=\!\!\!\!
\coprod_{d|(m,n)}\!\!\!\!
\coprod^{d}\!\!\!\!\!\!\!\!\!\!\!\!\!\!\!\!\!\!
\coprod_{\substack{ L\subset\Gamma \\
\ \ \ \ \ \ \ \ |\Gamma/L|=(mn)/d^2 }}
\!\!\!\!\!\!\!\!\!\!\!\!\!\!\!\!\!\!\!\!
\mathcal{F}(T_{d\cdot L})\big/\bigl(\Gamma/(d\cdot L)\bigr).
\end{equation*}
On the other hand, 
\begin{align*}
\bigl(\mathbb{T}\left(\frac{mn}{d^2}\right)(\mathbb{R}(d)\mathcal{F})\bigr)(T)
&=\Bigl[\!\!\!\!\!\!\!\!\!\!\!\!\!\!\!\!\!\!
\coprod_{\substack{ L\subset\Gamma \\\ \ \ \ \ \ |\Gamma/L|=(mn)/d^2 }}
\!\!\!\!\!\!\!\!\!\!\!\!\!\!\!\!\!\!
(\mathbb{R}(d)\mathcal{F})(T_L)/(\Gamma/L)\Big] \\
&=\!\!\!\!\!\!\!\!\!\!\!\!\!\!\!\!\!\!
\coprod_{\substack{ L\subset\Gamma \\\ \ \ \ \ \ |\Gamma/L|=(mn)/d^2}}
\!\!\!\!\!\!\!\!\!\!\!\!\!\!\!\!\!\!
\bigl[\mathcal{F}(T_{d\cdot L})\big/\bigl(L/(d\cdot L)\bigr)\bigr]
\big/(\Gamma/L) \\
&=\!\!\!\!\!\!\!\!\!\!\!\!\!\!\!\!\!\!
\coprod_{\substack{ L\subset\Gamma \\\ \ \ \ \ \ |\Gamma/L|=(mn)/d^2 }}
\!\!\!\!\!\!\!\!\!\!\!\!\!\!\!\!\!\!
\mathcal{F}(T_{d\cdot L})\big/\bigl(\Gamma/(d\cdot L)\bigr).
\end{align*}
Thus combining the above calculations, we have our formula (4-12). 
\end{proof}
Theorem B is a special case of Theorem 4.4 when
$\mathcal{F}(\Sigma)=\text{Map}_{\text{orb}}(\Sigma, M/G)$.

By a general procedure, Theorem 4.4 implies the following ``formal''
Euler product of operators:
\begin{equation}
\sum_{n\ge1}\frac{\mathbb{T}(n)}{n^s}
=\!\!\!\!\!
\prod_{p:\text{ prime}}
\!\!\!\!\!
\bigl(1-\mathbb{T}(p)p^{-s}
+p\cdot\mathbb{R}(p)p^{-2s}\bigr)^{-1},
\label{4.14}
\end{equation}
on functors $\mathcal{F}$.  However, the implications of this Euler product
formula in our present context is not clear.


\begin{thebibliography}{DHVW}

\bibitem{BL}  L.~Borisov and A.~Libgober, 
\emph{Elliptic genera of singular varieties}, 
Duke Math. J., 116 (2003), 319--351. 

\bibitem{DHVW} L.~Dixon, J.~Harvey, C.~Vafa and E.~Witten,
\emph{Strings on orbifolds}, Nuclear Physics, B 261, (1985), 678--686. 

\bibitem{DMVV} R.~Dijkgraaf, G.~Moore, E.~Verlinde, and H.~Verlinde,
\emph{Elliptic genera of symmetric products and second quantized
strings}, Comm. Math. Phys. 185 (1997), 197--209. 

\bibitem{Ga} N.~Ganter, \emph{Orbifold genera, product formulas and power operations}, Adv. math. 205 (2006), no. 1, 84--133. arXiv:math/0407021. 

\bibitem{GLM} S.~M.~Gusein-Zade, I.~Luengo, and A.~Melle-Hern\'andez,
\emph{On the power structure over the Grothendieck ring of varieties and its applications}, arXiv:math/0605467. 

\bibitem{H} F.~Hirzebruch, \emph{Elliptic genera of level $N$ for complex manifolds}, Diff. Geom. Methodsin Theoretical Physica, Kluwer Dordrecht, 1988, pp.37--63. 

\bibitem{HH} F.~Hirzebruch and H.~H\"ofer, 
\emph{On the Euler number of an orbifold}, Math.~ Annalen, 
286, (1990) 255--260

\bibitem{Land} P.S.~Landweber (ed.), \emph{Elliptic Curves and Modular Forms in Algebraic Topology}, Lecture Notes in Math. (Proceedings, Princeton 1986), vol 1326, Springer-Verlag, New York, 1988. 

\bibitem{L} S.~Lang, \emph{Introduction to Modular Forms}
Grund. Math. Wiss. 222, Springer-Verlag, New York, 1976. 

\bibitem{May} J.P.~May, \emph{Equivariant homotopy and cohomology theory}, CBMS regional Conferenbce Series in Mathematics, 91, Amer. math. Soc., Providence, RI, 1996.

\bibitem{M} I.~Moerdijk, \emph{Orbifolds as Groupoids: an Introduction}, Orbifolds in mathematics and physics (Madison WI, 2001), Contemp. Math. vol. 310, Amer. Math. Soc., Providence, RI, 2002, pp.205--222. 

\bibitem{T1} H.~Tamanoi, \emph{Generalized orbifold Euler 
characteristic of symmetric products and equivariant Morava K-theory}
Algebraic and Geometric Topology, 1 (2001) 115--141. 

\bibitem{T2} H.~Tamanoi, \emph{Generalized orbifold Euler 
characteristics of symmetric orbifolds and covering spaces}
Alg. Geom. Topology, 3 (2003), 791--856. 

\bibitem{T3} H.~ Tamanoi, \emph{Elliptic Genra and Vertex Operator Super-Algebras}, Lecture Notes in Math. vol 1704, Springer-Verlag, 1999. 

\bibitem{W} W.~Wang, \emph{Equivariant K-theory, wreath products, 
and Heisenberg algebra}, Duke Math. J. 103 (2000), 1--23. 

\bibitem{WZ} W.~Wang and J.~Zhou, \emph{Orbifold Hodge numbers of the wreath product orbifolds}, J. Geom. Phys. 38(2001) 153--170. arXiv:math/0005124.

\bibitem{Wi} E.~Witten, \emph{The index of Dirac operators in loop space},
Lecture Notes in Math. (Proceedings, Princeton, 1986, P.S.Landweber, ed.), vol 1326, Springer-Verlag, New York, 1988, pp 161--181. 

\end{thebibliography}
\end{document}